\documentclass[10pt,draftcls,onecolumn]{IEEEtran}

\IEEEoverridecommandlockouts

	\usepackage{graphicx}

\makeatletter
\def\maxwidth{
  \ifdim\Gin@nat@width>\linewidth
    \linewidth
  \else
    \Gin@nat@width
  \fi
}
\makeatother

\usepackage{placeins}
\usepackage{amsmath}
\usepackage{amsfonts}
\usepackage{amsthm}
\usepackage{amssymb}
\usepackage{multirow}
\usepackage{floatrow}
\usepackage[noadjust]{cite}
\usepackage{fixltx2e}
\usepackage{verbatim}
\usepackage{color}
\newfloatcommand{capbtabbox}{table}[]

\usepackage{blindtext}

\newtheorem{theorem}{Theorem}

\newtheorem{corollary}{Corollary}
\newtheorem{definition}{Definition}
\newtheorem{lemma}{Lemma}

\newtheorem{remark}{Remark}

\begin{document}

\title{Joint Centrality Distinguishes Optimal Leaders in Noisy Networks}
\author{Katherine Fitch and  Naomi Ehrich Leonard
\thanks{Supported in part by  NSF Graduate Research Fellowship DGE 1148900, NSF grant ECCS-1135724, ONR grant N00014-09-1-1074, and ARO grant W911NG-11-1-0385.
}
 \thanks{K. Fitch and N. E. Leonard are with the Department of Mechanical and Aerospace Engineering, Princeton University,
Princeton, NJ 08544, USA.\{ kfitch, naomi\}@princeton.edu}%
}
\maketitle

\begin{abstract}
We study the performance of a network of agents tasked with tracking an external unknown signal in the presence of stochastic disturbances and under the condition that only a limited subset of agents, known as leaders, can measure the signal directly.  We investigate the optimal leader selection problem for a prescribed maximum number of leaders, where the optimal leader set minimizes total system error defined as  steady-state variance about the external signal. In contrast to previously established greedy algorithms for optimal leader selection, our results rely on an expression of total system error in terms of properties of the underlying network graph.  We demonstrate that the performance of any given set of \textcolor{black}{noise-free} leaders depends on their influence as determined by a new graph measure of centrality of a set.  We define the {\em joint centrality} of a set of nodes in a network graph such that a \textcolor{black}{noise-free} leader set with maximal joint centrality is an optimal leader set. 
In the case of a single leader, we prove that the  optimal leader is the node with maximal information centrality \textcolor{black}{for both the nois-corrupted and noise-free leader cases}.  In the case of multiple leaders, we show that the nodes in the optimal \textcolor{black}{noise-free} leader set balance high information centrality with a coverage of the graph.  For special cases of graphs, we solve explicitly for optimal leader sets.  We illustrate  with examples. 
\end{abstract}

\section{Introduction}

Analysis of networked multi-agent system dynamics has generated substantial research interest in recent years \cite{OlfatiSaber2004,ren2005,Xiao2007}.  This is largely due to the broad range of applications for which the theory can be applied, including, for example, design of vehicle networks \cite{ren}, analysis of social networks \cite{Jadbabaie2012}, investigation of collective animal behavior \cite{Starlings2013}, and more. Often in these applications the network must learn   an external signal, for example, in the case of a sensor network using consensus to estimate an environmental signal \cite{olfatisaber2005}. However, when the signal is costly to sample, e.g. because of energy consumption costs or costs associated to acquiring the necessary sensory or processing capability, it may become impractical for all agents in the network to measure the signal directly.  If inter-agent sensing or communication is relatively inexpensive, then a more efficient solution involves a limited subset of agents, called leaders, measuring the signal directly, with   
 the remaining agents, called followers, learning the signal through network connections.   
 In this paper we address the problem of selecting leaders, as a function of the network graph, to maximize network accuracy in tracking the external signal.
 
The problem is motivated by the design of high performing engineered networks such as sensor networks, as well as by finding conditions under which biological networks such as animal groups perform highly.   For example, in migratory herds, the animals must learn, agree on, and move together along a single migration route.   
It is likely that only a  subset of animals invest in a direct measurement of the route, particularly when it is easier to rely on observations of neighbors \cite{Guttal2010}.  
In \cite{Guttal2010} it was shown that the emergence of leaders and followers within a large, mobile population is an evolutionarily stable solution for sufficiently high investment cost of sampling the migration route.   In \cite{Pais2012}, the authors used a mathematical model to analyze this evolutionary dynamic and to compute the location of emergent leaders as a function of the network graph and the investment cost.  The model yields a distributed adaptive dynamic for taking on leadership in this context; however, the evolutionary dynamics do not guarantee a steady-state solution that is optimal for the herd.  



In the present paper, we study the leader-follower network  dynamics subject to stochastic perturbations (\cite{WangTan2009,Patterson2010}), examining cases in which \textcolor{black}{there are one and two noise-corrupted leaders and in which there are any number of noise-free leader nodes}.  Our objective is to make rigorous how a leader set, as a function of properties of the network graph, affects the total system error of the group defined as the steady-state variance of the system about the external signal. Total system error can also be viewed as a measure of coherence, equivalently the $H_2$ norm of the system dynamics \cite{Patterson2010,young2010}.

To this end we develop a means of quantifying the combined influence of a set of leader nodes in a network on the total system error in the leader-follower dynamic. Intuitively,  this influence should correspond to some notion of  centrality of a set of nodes since a leader set that gives low system error must be well connected to other nodes in the network. Different types of centrality of a set of nodes were defined in \cite{groupcent}, where the authors quantified degree, closeness, betweenness and flow centralities of sets of nodes by extensions of the  definitions for individuals.   Illustrative examples were  used in \cite{groupcent} to explore the relationship between those measures and network properties. In contrast to  the literature, we derive a measure of centrality of a set of nodes,  called {\em joint centrality}, by examining the performance measure, i.e., total system error, and expressing performance in terms of graph measures. 

We note that applications of measures of centrality of a set of nodes include a broad range of research areas from emergency response management \cite{emergresp} to a network connectivity analysis of the quality of innovative ideas \cite{innovation}. 

Much of the recent research related to leader-follower multi-agent systems with stochastic dynamics has been focused on the development of off-line leader selection algorithms that seek to find the leader set that minimizes total system error \cite{Patterson2010,FarLin1,LinFar2,SuperOpt,Jovanovic2013}. These algorithms have been designed to be computationally efficient in approximating optimal solutions with proven bounds on the total system error relative to the optimal value of error. Many of the algorithms are iterative, adding to the leader set one agent at a time. This approach may preclude finding the optimal solution  since the optimal set of $l$ leaders does not necessarily include the optimal set of $m$ leaders, $l>m$. The authors of \cite{Jovanovic2013} address this issue  by considering a ``swap'' step within their iterative algorithm.   

Our contributions are fourfold.  First, we provide a new approach to solving the optimal leader selection problem in terms of network graph measures.   In general, our approach reduces computational complexity significantly as compared to the brute force computation.       Second, we define a new notion of centrality of a set of nodes in an undirected,  connected graph, that we call joint centrality.   
For the leader-follower network tracking dynamics, we show that the total system error  is proportional to the joint centrality of the  leader set when the leaders are noise-free.  The joint centrality of a set of nodes depends on the information centrality of the nodes and the resistance and biharmonics distances between pairs of nodes in the set.   We show how to calculate joint centrality using  entries from  submatrices of the pseudoinverses of the Laplacian and squared Laplacian.   Third, we consider the case of noise-corrupted leaders and we derive a modified notion of joint centrality, showing, in the cases of one and two noise-corrupted leaders, that total system error is proportional to the modified joint centrality of the leaders.   Fourth, we prove the explicit solution to the optimal leader selection problem in the case of cycle graphs and path graphs.   Further, in the case of one noise-free or noise-corrupted leader, we prove that the optimal leader is the agent with  maximal information centrality.      A preliminary version of results in the paper, for the cases of one and two leaders, appears in \cite{Fitch2013}.

The paper is organized as follows. In Section~\ref{sec:model}, we introduce the network model dynamics and define the optimal leader selection problem. We  review information centrality, resistance distance, biharmonic distance and other  properties of the Laplacian in Section~\ref{sec:graphs}.   In Section~\ref{general} we derive total system error for the general case of $m$ noise-free leaders,   define and interpret joint centrality of $m$ nodes and provide an illustrative example.   In Section~\ref{sec:interp} we \textcolor{black}{ provide an interpretation of joint centrality. In addition} we consider both the noise-free and noise-corrupted cases for two leaders and derive exact solutions to special cases, prove the relationship between an optimal single leader and information centrality, \textcolor{black}{and explain the connection to the problem of controllability on networks}.   We show an example application of joint centrality in Section~\ref{sec:special}.   We conclude with a discussion in Section~\ref{sec:final}.

\section{Model and Problem Statement}
\label{sec:model}

We consider a network of $n$ agents tasked with tracking an external signal from the environment.   We denote the external signal by $\mu \in \mathbb{R}$ and suppose it to be a constant.   Generalizations to vector-valued environmental signals are expected to be relatively straightforward and extensions to time-varying environmental signals are the topic of future work.

The state of agent $i$, for $i = 1, \ldots, n$, is  $x_i \in \mathbb{R}$, and it represents agent $i$'s estimate of the signal $\mu$.   
The state of the network is given by $\mathbf{x} = [x_1,x_2,...,x_n] \in \mathbb{R}^n$.     
Agent $i$ can measure the evolving relative state  $x_j - x_i$  for each agent $j$ in its set of neighbors $\mathcal{N}_i$.   The availability of these measurements to agent $i$ is the result of agent $i$ directly sensing the relative state of its neighbors, e.g., in the case that the state refers to position, or of neighbors  communicating the value of their state to agent $i$.    

The graph $\mathcal{G = (V,E},A)$ encodes the network topology.   Each agent corresponds to a node in the set $\mathcal{V} = \{1,2,...n\}$, and we will use the terms agents and nodes interchangeably.   $\mathcal{E \subseteq V \times V}$ is the set of edges, where  the edge $(i,j) \in \mathcal{E}$ if $j \in \mathcal{N}_i$.   The adjacency matrix is given by $A \in \mathbb{R}^{n \times n}$ where matrix element $a_{i,j}$ corresponds to the weight on edge $(i,j)$. 

We consider undirected, connected graphs.   Recent results on effective resistance in directed graphs  \cite{dres1,dres2} suggest the means to extend our theory in future work to the case of directed graphs. The undirected graph contains edge $(i,j)$, then $a_{i,j} = a_{j,i}  > 0$ and otherwise $a_{i,j}=0$.   The degree matrix $D$ is a diagonal matrix with entries $d_{i,i} = \sum_{j=1}^n a_{i,j}$. The  associated Laplacian matrix is defined as $L = D-A$.  

An agent $l \in \mathcal{V}$ is called a {\em leader} if it directly measures the external signal.  Let  $k_l > 0$ be the weight  that agent $l$ puts on its signal measurement.  Any agent that is not a leader is called a {\em follower}.  Let the set of leaders be denoted $S$ with cardinality $m$ and the set of follower nodes, denoted by $F$, be the complement of $S$ with cardinality $n-m$. Summation over $s$ denotes summation over the leader set, while summation over $i$ denotes summation over the entire set of leaders and followers. We use the index $l_1$ when it is necessary to identify one leader apart from the rest of the leader set. 

Throughout the paper, when a set $S$ of $m$ nodes is identified, we will assume they are the first $m$ nodes in an ordering of the $n$ nodes.  Accordingly, we will denote the partition of an $n \times n$ matrix $B$ as
\begin{align}
B = \begin{bmatrix}B_{S} && B_{SF}\\ B_{FS} && B_{F} \end{bmatrix}, \label{mpartition} 
\end{align}
where  $B_S$ is an $m\times m$ matrix corresponding to nodes in set $S$, and $B_F$ is an $(n-m)\times(n-m)$ matrix corresponding to the remaining nodes.   We will further let $l_1$ be the first node in the ordered set $S$.
We will denote the Moore Penrose pseudoinverse of a matrix $B$  by $B^+$ and the conjugate transpose of $B$ by $B^*$.   We let $\mathbf{1_n}$ be the vector of $n$ ones and $\mathbf{e}_j$ be the standard basis vectors for $\mathbb{R}^n$.

We  assume  that all leaders apply the same weight $k$ to their measurement of the external signal, i.e.,  $k_i = k > 0$ for $i \in S$ and $k_i = 0$ for $i \in F$.    We assume that stochastic disturbances enter the dynamics as additive noise.   We model the dynamics for each agent $i \in \mathcal{V}$ by the following stochastic process:
\begin{align}
dx_i =  -k_i(x_i -\mu) dt -L_i\mathbf{x} dt+ \sigma dW_i, \label{dyn1} 
\end{align}
where $L_i$ is the $i$th row of the Laplacian $L$, and $\sigma dW_i$  represents increments drawn from independent Wiener processes with standard deviation $\sigma$.   

In the case that $k < \infty$, the dynamics of the leaders and followers are all noise corrupted.
In \cite{Jovanovic2013}, it was demonstrated that in the limit as $k \rightarrow \infty$, i.e., in the case that leaders apply an arbitrarily large weight to tracking the external signal, the dynamics (\ref{dyn1}) describe the case of noise-free leaders.    Thus, our model (\ref{dyn1}) describes  both  cases of noise-corrupted leaders ($k<\infty$) and noise-free leaders ($k \rightarrow \infty$).

To write (\ref{dyn1}) in vector form let $K\in \mathbb{R}^n$ be the diagonal matrix with elements $k_i$, let $M=L+K$ and without loss of generality let $\mu = 0$. Then  (\ref{dyn1}) becomes
\begin{align}
d\mathbf{x} = -M \mathbf{x} dt+ \sigma d\mathbf{W}. \label{vec1} 
\end{align}
Since we have assumed that $\mathcal{G}$ is connected, $-M$ is Hurwitz so long as $k_i=k>0$ for some agent $i$, i.e., $S$ is nonempty. 

Thus, for nonempty $S$,  $\mathbf{x}$ will converge to a steady-state  distribution about the value of the external signal, and the steady-state covariance matrix $\Sigma$ of $\mathbf{x}$ is the solution to the Lyapunov equation
\begin{align}
M\Sigma + \Sigma M^T = \sigma^2 I. 
\end{align}
The steady-state variance of  $x_i$ is $\Sigma_{i,i}$, the corresponding diagonal element of $\Sigma$. Since the external signal is assumed to be constant, the system will converge to a steady-state distribution about the value of the external signal even if the nodes chosen as leaders do not guarantee system controllability.  

Following \cite{Patterson2010,SuperOpt}, we define {\em total system error} as $\text{tr}(\Sigma) = \sum_{i=1}^n \Sigma_{i,i}$.   We define group performance as the inverse of total system error, which measures {\em network tracking accuracy}.

By \cite{Arnold} we have that the covariance matrix of (\ref{vec1}) is 
\begin{align}
\text{Cov}(\mathbf{x}(t),\mathbf{x}(t)) = \sigma^2 \int_0^t e^{-M(t-\tau)} e^{-M^T(t-\tau)} d\tau. \label{cov1} 
\end{align}
Given that $\mathcal{G}$ is undirected, the Laplacian matrix $L$ will be symmetric and it follows that $M$ will be symmetric and normal.  Let the eigenvalues of $M$ be $\lambda_i$, $i \in \mathcal{V}$ with corresponding eigenvectors $\boldsymbol \nu_i$.   Let $\Lambda$ be the diagonal matrix with entries $\Lambda_{i,i} = \lambda_i$. Then there exists a unitary matrix $U$ such that $U^*MU=\Lambda$ and  (\ref{cov1}) can be written as
\begin{align}
\text{Cov}(\mathbf{x}(t),\mathbf{x}(t)) = \sigma^2(UR(t)U^*),
\end{align}
with
\begin{align}
R(t) := \int_0^t e^{-(\Lambda+\bar{\Lambda})(t-\tau)}d\tau.
\end{align}

From \cite{Poulakakis2010}, this gives
\begin{align}
[\text{Cov}(\mathbf{x}(t),\mathbf{x}(t))]_{i,j}=\sigma^2 \sum_{p=1}^n\frac{1-e^{-2\text{Re}(\lambda_p)t}}{2\text{Re}(\lambda_p)} \nu_i^{(p)}\bar{\nu}_j^{(p)}. 
\end{align}
Since $M$ is symmetric, all eigenvalues of $M$ will be real, and the steady-state variance of each node can be written as 
\begin{align}
\text{Var}(x_i)_{ss} = \Sigma_{i,i} = \sigma^2 \sum_{p=1}^n\frac{1}{2\lambda_p} |\nu_i^{(p)}|^2. \label{vars1}
\end{align}
Total system error follows from summing (\ref{vars1})  over all $i$,
\begin{align}
\sum_{i=1}^n \Sigma_{i,i} =  \sigma^2 \sum_{i = 1}^n \frac{1}{2\lambda_i}= \frac{\sigma^2}{2}\sum_{i = 1}^n M_{i,i}^{-1}. \label{arn} 
\end{align}
Total system error defines the coherence of the network, and is equivalent to the $H_2$ norm of the system \cite{Patterson2010,young2010}.   

We define the {\em optimal leader selection problem} as follows. 

\begin{definition}[Optimal leader selection problem]   Given $m$ and undirected, connected graph ${\cal G}$, find a set of $m$ leaders $S^*$ over all possible  sets $S$ of $m$ leaders that minimizes the total system error  (\ref{arn}) for the leader-follower network tracking dynamics (\ref{vec1}), i.e., find
\begin{align}
S^* = \arg \min_{S} \sigma^2 \sum_{i=1}^n \frac{1}{2\lambda_i}  = \arg \min_{S} \frac{\sigma^2}{2} \sum_{i=1}^n M_{i,i}^{-1}. \label{eq2}
\end{align}
\end{definition}

\section{Review of Properties of the Laplacian and Graph Theoretic Measures}
\label{sec:graphs}

We briefly review relevant graph theoretic measures and identities that will be applied in later sections.
We start with the notion of information centrality, which was first introduced by Stephenson and Zelen in \cite{infocent}. Information centrality can be understood by first defining the information in a path between any two nodes in $\mathcal{G}$ to be the inverse of the path length between those two nodes. Thus, the longer the path the less information in that path. Total information between nodes $i$ and $j$, denoted $I^{tot}_{i,j}$, is  the sum of the information in all paths connecting nodes $i$ and $j$. It was shown in \cite{infocent} that total information can be calculated without path enumeration by using  the group inverse of the Laplacian, which here is the pseudoinverse $L^+$:
\begin{align}
I_{i,j}^{\text{tot}} = (L^+_{i,i} + L^+_{j,j} -2L^+_{i,j})^{-1}.
\end{align}
Information centrality for node $i$, denoted $c_i$,  is defined as the harmonic average of total information between node $i$ and all other nodes in $\mathcal{G}$:
\begin{align}
c_i = \left(\frac{1}{n} \sum_{j=1}^n \frac{1}{I_{i,j}^{\text{tot}}} \right)^{-1}. \label{infocen}
\end{align}

In \cite{NodeCert}, Poulakakis et al.\  evaluated the certainty of each node $i$ in a network of decision-makers accumulating stochastic evidence towards a decision.  This certainty, denoted $\mu_i$, is defined as the inverse of the difference between the variance of the state $x_i$ about the reference signal and the minimum achievable variance as $t \rightarrow \infty$. The authors apply the notion of information centrality to directly interpret $\mu_i$ in terms of structural properties of the underlying communication graph.  It was proven that 
\begin{align}
\frac{1}{\mu_i} = \frac{\sigma^2}{2} L_{i,i}^+ = \frac{\sigma^2}{2} \left(\frac{1}{c_i}-\frac{K_f}{n^2} \right), \label{Lp}
\end{align}
where $K_f$ is the Kirchhoff index of $\mathcal{G}$. The identity (\ref{Lp}) implies that the ordering of nodes by certainty is equal to the ordering of nodes by information centrality. We  show in later sections that information centrality also plays a critical role in the solution to the optimal leader selection problem.

The total information between any two nodes $i$ and $j$ is closely related to the resistance distance between them, denoted $r_{i,j}$.  Resistance distance between nodes in the undirected graph $\mathcal{G}$ is defined as the resistance distance between the corresponding two nodes in the electrical network analog to the graph $\mathcal{G}$.   By \cite{resdist} for an undirected graph $\mathcal{G}$ 
\begin{align} r_{i,j} & = L_{i,i}^+ +L_{j,j}^+ - 2L_{i,j}^+ = {I_{i,j}^{\text{tot}}}^{-1}.\label{res} \end{align}
It follows that  
\begin{align}
\sum_{i=1}^n r_{i,j} = \frac{n}{c_j}. \label{rescen}
\end{align}

An additional measure with similar form to that of resistance distance is the recently derived notion of biharmonic distance, $d_B$ \cite{biharmonic}.   This measure has been used to quantify distance between two points $v_i, v_j$ on the surface of a discrete 3D  mesh: 
\begin{align}
d_B(v_i,v_j)^2 = g_d(i,i) + g_d(j,j) -2g_d(i,j), \label{bd}
\end{align}
where $g_d$ is the discrete Green's function of the discretized, normalized bilaplacian $\tilde L^2$, equivalent to the pseudoinverse of $\tilde L^2$, and $\tilde{L}$ is the normalization of Laplacian $L$.  We define the {\em biharmonic distance between two nodes $i$ and $j$ in the graph ${\cal G}$}, which we denote $\gamma_{i,j}$, analogously   without normalizing $L$:
 \begin{align}
\gamma_{i,j}& = L^{2+}_{i,i} + L^{2+}_{j,j} - 2L^{2+}_{i,j} =\sum_{l=1}^n (L_{l,i}^+ - L_{l,j}^+)^2 \nonumber \\
&=  (\mathbf{e}_i-\mathbf{e}_j)^T L^{2+} (\mathbf{e}_i-\mathbf{e}_j). \label{s1}
\end{align} 

We observe that the biharmonic distance $\gamma_{i,j}$ of (\ref{s1}) is very similar to resistance distance $r_{i,j}$ of (\ref{res}) with the difference being the use of the pseudoinverse of $L^2$ in the definition of $\gamma_{i,j}$ as compared to the pseudoinverse of $L$ in the definition of $r_{i,j}$. Since $L^{2}$ is symmetric and positive semi-definite, we immediately have that $\gamma^{1/2}$ is a metric. In fact, it can be viewed as a Manahalobis distance, which in this case describes a dissimilarity measure between two vectors from a single distribution with  covariance matrix $L^2$.  Let  $\Gamma$ be the matrix with elements $\gamma_{i,j}$. 

For completion, we note that both resistance distance and biharmonic distance between nodes can be written in terms of the eigenvalues $\lambda_i$ and eigenvectors $\boldsymbol{\nu}_i$ of the Laplacian   $L$:
\begin{align}
r_{i,j} = \sum_{l=2}^n \frac{1}{\lambda_l}(\nu_l^i - \nu_l^j)^2, \\
\gamma_{i,j} = \sum_{l=2}^n \frac{1}{\lambda_l^2}(\nu_l^i - \nu_l^j)^2.
\end{align}

Finally, the following properties of $L^+$ will be applied in proofs (see \cite{NodeCert} for details):
\begin{align}
&LL^+ = L^+L = I_n - \frac{1}{n} \mathbf{1_n 1_n}^T, \label{Lp1 }\\ 
&\mathbf{1_n}^TL^+ = L^+ \mathbf{1_n }= 0, \label{Lp2} \\   
&\mathrm{Tr}(L^+) = \frac{K_f}{n}. \label{Lp3}
\end{align}

\section{Joint centrality and the optimal $m$ noise-free leaders} \label{general}
In this section, we prove our main result on the general solution of the optimal leader selection problem by deriving an explicit expression for total system error with $m$ noise-free leaders in terms of properties of the underlying graph.   Before stating the theorem, we first define the {\em joint centrality of a set of $m$ nodes} in a network graph. 

\begin{definition}[Joint centrality]  \label{jc}  Let $\mathcal{G}$ be an undirected, connected graph of order $n$.   Given integer $m<n$, let $S$ be the set of any $m$ nodes in ${\cal G}$.    Choose an arbitrary element $l_1 \in S$.  Let $N$ be an $n \times n$ matrix with elements of  $N^{-1}$ given by
\begin{align}
N^{-1}_{i,j} = L^+_{i,j}-L^+_{i,l_1} - L^+_{j,l_1}+L^+_{l_1,l_1}. \label{nterms}
\end{align}
Following (\ref{mpartition}),  $N^{-1}_{S \backslash l_1}$ is the $(m-1)\times (m-1)$ submatrix of $N^{-1}$ corresponding to the elements of $S$ less the first element $l_1$.  Let $G = \left(N^{-1}_{S \backslash l_1}\right)^{-1}$ and  $\bar G  = \begin{bmatrix} 0 & 0\\0 & G \end{bmatrix} \in \mathbb{R}^{m \times m}$.  Let $Q = \bar G \Gamma_S$, where $\Gamma$ is given by (\ref{s1}).   
The {\em joint centrality of set $S$ in ${\cal G}$} is defined as
\begin{align}
\rho_S =  n\Big(&\frac{K_f}{n} +n \det( G) \det( L_S^+)  
+ \frac{1}{2} \mathrm{Tr}(Q) - \mathbf{1}_n^TQ\mathbf{e}_{l_1}\Big)^{-1}. \label{rhom}
\end{align}
\end{definition}
 
\begin{theorem} \label{thm1} Let $\mathcal{G}$ be an undirected,  connected graph of order $n$.  Let  $S$ be a set of $m$ noise-free leaders. Then, the total system error (\ref{arn}) for the system dynamics (\ref{vec1}) is
\begin{align}
\sum_{i=1}^n \Sigma_{i,i} = \frac{\sigma^2}{2}\Big( \frac{n}{\rho_S}\Big), \label{first}
\end{align}
where $\rho_S$ is the joint centrality of leader set $S$ given by (\ref{rhom}).
The optimal leader set is $S^* = \arg\max_{S} \rho_{S}$, the set of leader nodes with the maximal joint centrality.
\end{theorem}

We recall three lemmas that will be used in the proof of Theorem~\ref{thm1}.

\begin{lemma} \label{one_up} \cite{GenInvses} Let $\mathbf{z},\mathbf{y} \in \mathbb{R}^n$.   A rank-1 update $\mathbf{zy}^T$  for the Moore-Penrose pseudoinverse of a real valued matrix, $F\in \mathbb{R}^{n\times n}$, is given by
\begin{align} 
(F + \mathbf{zy}^T)^+ = F^+ +H \label{m_inv}
\end{align}
where 
\begin{align}
H = -\frac{1}{\|\mathbf{w}\|^2}\mathbf{vw}^T-\frac{1}{\|\mathbf{m}\|^2}\mathbf{mh}^T+ \frac{\beta}{\|\mathbf{m}\|^2\|\mathbf{w}\|^2}\mathbf{mw}^T \label{eqG}
\end{align}
and $\beta = 1+\mathbf{y}^TF^+\mathbf{z}$, $\mathbf{v} = F^+\mathbf{z}$, $\mathbf{h} = (F^+)^T\mathbf{y}$, $\mathbf{w} = (I - FF^+)\mathbf{z}$, and $\mathbf{m} = (I-F^+F)^T \mathbf{y}$.
\end{lemma}

\begin{lemma} \label{gen_up}\cite{wb1950} Let $X \in \mathbb{R}^{n \times n}$, $Z \in \mathbb{R}^{m \times m}$, $U \in \mathbb{R}^{n \times m}$ and $V \in \mathbb{R}^{m \times n}$ such that  $X$, $Z$ and $X+UZV$ are nonsingular.   Then, $(X+UZV)^{-1}$ can be written as
\begin{align}
(X+&UVZ)^{-1} = X^{-1}-X^{-1}U(Z^{-1}+VX^{-1}U)^{-1}VX^{-1}. \label{inv3}
\end{align}
\end{lemma}

\begin{lemma} \label{det_partition} \cite{detpart} The determinant of a bordered matrix can be computed as follows
\begin{align}
\left| \begin{array}{cc}
X & \mathbf{u} \\
\mathbf{v}^T & d  \end{array} \right| = d |X| -\mathbf{v}^T(\mathrm{adj }X)\mathbf{u},
\end{align}
where $X \in \mathbb{R}^{p \times p}$, $\mathbf{u}, \mathbf{v} \in \mathbb{R}^{p}$, and $d \in \mathbb{R}$.
 \end{lemma}
\proof 
(Theorem~\ref{thm1}).  We begin by studying terms in the total system error for finite $k>0$ and then evaluate in the limit as $k \rightarrow \infty$. From (\ref{arn}), the total system error is proportional to  Tr$(M^{-1})$ where $M = L+K$. Let  $K_1$ be the diagonal matrix with $k$ in the first diagonal element and zeros elsewhere and let $K_{m-1} = K - K_1$. 
We derive an expression for Tr$(M^{-1})$ by calculating two successive updates to $L^+$. We first show that if we define $N = L + K_1$, and thus $M = N + K_{m-1}$, then $N^{-1}$ satisfies (\ref{nterms})  for $k \rightarrow \infty$. 

Let $\mathbf{e} = \mathbf{d}$ be vectors of length $n$ with $\sqrt{k}$ in the $l_1$ (first) entry and zeros elsewhere \textcolor{black}{where $l_1$ is a member of the leader set. Note that the choice of $l_1$ will not effect the value of joint centrality for a given leader set.} Then $N^{-1} = (L + K_1)^{-1} = (L + \mathbf{e} \mathbf{d}^T)^{-1}$.   Applying Lemma~\ref{one_up} we get that $(L + \mathbf{e} \mathbf{d}^T)^{-1} = L^+ + H$, with $H$ given by  (\ref{eqG}) such that
\begin{align}
N^{-1} =  &L^{+}  -L_{l_1}^+\mathbf{1_n}^T -\mathbf{1_n}L_{l_1}^{+T} + \frac{ (1+kL_{l_1,l_1}^+)}{k}\mathbf{1_n}^T\mathbf{1_n}. \label{nmat}
\end{align}
Taking the limit as $k \rightarrow \infty$, the elements of $N^{-1}$ can be written as (\ref{nterms}).

Let $U = [-\sqrt{k}\mathbf{e}_2, \ldots, -\sqrt{k}\mathbf{e}_m] \in \mathbb{R}^{n \times (m-1)}$, let $V = U^T$ and let $\mathbb{I}_{m-1} \in \mathbb{R}^{(m-1)\times(m-1)}$ be the identity matrix.   Then, $M^{-1} = (N + K_{m-1})^{-1} = (N + U\mathbb{I}V)^{-1}$.
Applying Lemma \ref{gen_up} we get that
\begin{align}
(N+&U\mathbb{I}V)^{-1} = N^{-1} - N^{-1}U(\mathbb{I}+VN^{-1}U)^{-1}VN^{-1}.
\end{align}
 
Let $G = ({N^{-1}}_{S\backslash l_1})^{-1}$ as in Definition~\ref{jc}.  Then if we take the limit as $k \rightarrow \infty$, sum the diagonal elements of $M^{-1} = (N+UIV)^{-1} $, and apply the identities  (\ref{Lp2}) and (\ref{Lp3}) we get
\begin{align}
\sum_{j=1}^n M_{j,j}^{-1} 
=& \frac{K_f}{n} + nL_{l_1,l_1}^+ -    \!\! \sum_{s_1,s_2 \in S \backslash \{l_1\}} \sum_{i=1}^n G_{s_1,s_2} 
\Big( L_{l_1,l_1}^+(L_{l_1,l_1}^+- L_{l_1,s_1}^+ - L_{l_1,s_2}^+) +  L_{l_1,s_1}^+L_{l_1,s_2}^+ +
 \nonumber \\
 &  \frac{1}{2} \left[(L_{i,l_1}^+ - L_{i,s_1}^+)^2 + (L_{i,l_1}^+ - L_{i,s_2}^+ )^2 \! - (L_{i,s_1}^+ -L_{i,s_2}^+)^2\right]\Big). \label{gen1}
\end{align}

Consider the square bracketed terms of  (\ref{gen1}) in which we observe the emergence of biharmonic distance, $\gamma$. 
Substituting (\ref{s1}) and  defining $\bar{G}$ as in Definition~\ref{jc} we get
\begin{align}
&  \!\! \sum_{s_1,s_2 \in S \backslash \{l_1\}} \sum_{i=1}^n G_{s_1,s_2} 
  \; \frac{1}{2} \left[(L_{i,l_1}^+ - L_{i,s_1}^+)^2 + (L_{i,l_1}^+ - L_{i,s_2}^+ )^2 \! - (L_{i,s_1}^+ -L_{i,s_2}^+)^2\right]\Big) 
= -\frac{1}{2} \mathrm{Tr}(\bar{G}\Gamma_S) +\mathbf{1}_n^T[\bar{G} \Gamma_S]\mathbf{e}_{l_1}. \label{simp2}
\end{align}

Additional simplification is made by applying Lemma \ref{det_partition} to the middle terms on the right hand side of (\ref{gen1}).  We get
\begin{align}
n&L_{l_1,l_1}^+- n\sum_{s_1,s_2 \in S \backslash \{l_1\}} G_{s_1,s_2} 
\Big( L_{l_1,l_1}^+(L_{l_1,l_1}^+ - L_{l_1,s_1}^+ - L_{l_1,s_2}^+) +  L_{l_1,s_1}^+L_{l_1,s_2}^+ \Big) =\nonumber \\
& \frac{n}{\det(G^{-1})}\Big(L_{l_1,l_1}^+\det(G^{-1}) - \sum_{s_1,s_2 \in S \backslash \{l_1\}} C_{{N^{-1}}_{s_1,s_2}} \left[ L_{l_1,l_1}^+(L_{l_1,l_1}^+- L_{l_1,s_1}^+   - L_{l_1,s_2}^+) +  L_{l_1,s_1}^+L_{l_1,s_2}^+\right] \Big) \label{middle}
\end{align}
where $C_{N^{-1}}$ is the cofactor matrix of ${N^{-1}}_{S\backslash l_1} = G^{-1}$. We then let $\mathbf{L}_{l_1,s_i}^+ = [L_{l_1,s_1}^+,...,L_{l_1,s_m-1}^+]^T$ and $\mathbf{L}_{l_1,l_1}^+ = [L_{l_1,l_1}^+,...,L_{l_1,l_1}^+]^T$ to be  vectors in $\mathbb{R}^{m-1}$ and apply Lemma \ref{det_partition} to rewrite the expression in (\ref{middle}) as
\begin{align}
 \frac{n}{\det({N^{-1}}_{S\backslash l_1})}
\left| \begin{array}{cc}
\  {N^{-1}}_{S\backslash l_1}& \mathbf{L}_{l_1,l_1}^+-{\mathbf{L}_{l_1,s_i}^+} \\
\ {\mathbf{L}_{l_1,l_1}^+}-{\mathbf{L}_{l_1,s_i}^+}& L_{l_1,l_1}^+
 \end{array} \right|.  \label{detproof1}
\end{align}

Using (\ref{nterms}) we  expand the determinant in (\ref{detproof1}) and perform algebraic manipulation to show that (\ref{detproof1}) simplifies to 
 \begin{align}
n \det( G)\det(L_S^+). \label{simp1}
 \end{align}
 Thus, 
 \begin{align}
\sum_{i=1}^n \Sigma_{i,i} &= \frac{\sigma^2}{2}\Big(\frac{K_f}{n} +n \det( G) \det( L_S^+)+ \frac{1}{2} \mathrm{Tr}(\bar{G}\Gamma_S) 
- \mathbf{1}_n^T[\bar{G} \Gamma_S]\mathbf{e}_{l_1} \Big) = \frac{\sigma^2}{2}\Big(\frac{n}{\rho_S}\Big)
\end{align}
where $\rho_S$ is defined by (\ref{rhom}).
 \endproof
 
 \section{\textcolor{black}{Interpretation}}
 \label{sec:interp}
 
\textcolor{black}{ In this section we provide interpretation and intuition behind the measure of joint centrality. In addition we provide an illustrative example and insights from the special cases of one and two noise-corrupted and noise-free leaders. We begin with two remarks.}

\begin{remark}  Using Theorem~\ref{thm1} to compute the total system error in terms of joint centrality of the $m$ leader nodes provides a significant reduction in computation as compared to using the definition of total system error (\ref{arn}).  Using joint centrality one only needs to compute the inverse of two $n \times n$ matrices $L^+$ and $L^{2+}$ and then for each candidate set of leaders the inverse of an $(m-1) \times (m-1)$ matrix.  This is in contrast to using the definition (\ref{arn}), which requires computing the  inverse of the  $n\times n$ matrix $M$ for each candidate set of leaders.

\end{remark}

\begin{remark}  \label{rem2} Theorem~\ref{thm1} reveals how the solution to the optimal leader selection problem is an optimal trade-off between high information centrality of  the leader nodes and high resistance distances and biharmonic distances between leader nodes.    To see this we examine the terms in  (\ref{rhom}) for joint centrality $\rho_S$.

 First, the elements of $N^{-1}$ given by (\ref{nterms}) depend on resistance distances:
\[
N^{-1}_{i,j} = \frac{1}{2} ( r_{i,l_1} + r_{j,l_1} - r_{i,j} ).
\]
Thus $N^{-1}_{i,j}$ quantifies a joint resistance distance between  a pair of nodes $i$, $j$ and $l_1$,  Then, $\det(G) =  (\det(N^{-1}_{S \backslash l_1}))^{-1}$ depends on these joint resistance distances among leaders.   

Second,  by (\ref{Lp})  each diagonal element of $L^+_S$ corresponds to a leader node and depends directly on the inverse of its information centrality as follows:
\[
L^+_{s,s} = \frac{1}{c_s} - \frac{K_f}{n^2}.
\]
By (\ref{res}) the off diagonal elements of $L^+_S$ depend on information centralities and resistance distances between leaders:
\[
L_{s,t}^+ = \frac{1}{2} \left(\frac{1}{c_s} + \frac{1}{c_t} - r_{s,t} - 2\frac{K_f}{n^2}  \right).
\]
Maximizing $\rho_S$ requires a small $\mathrm{det}(G)\mathrm{det}(L^+_S)$, which suggests a key trade-off between  high information centrality of leaders and high resistance distances between leaders.


The term $\mathrm{Tr}(Q)$ in (\ref{rhom}) is the sum of products of the biharmonic distance between pairs of leader nodes (from $\Gamma_S$), and terms in $G$.  
Since $\mathrm{Tr}(Q)$ is negative, maximizing joint centrality requires high biharmonic distance between pairs of leader nodes.   Thus, the optimal leader set trades off high information centrality with a coverage of the graph as made rigorous by the joint centrality measure.




\end{remark}

\begin{figure}[ht!]
\centering
\includegraphics[width=3in]{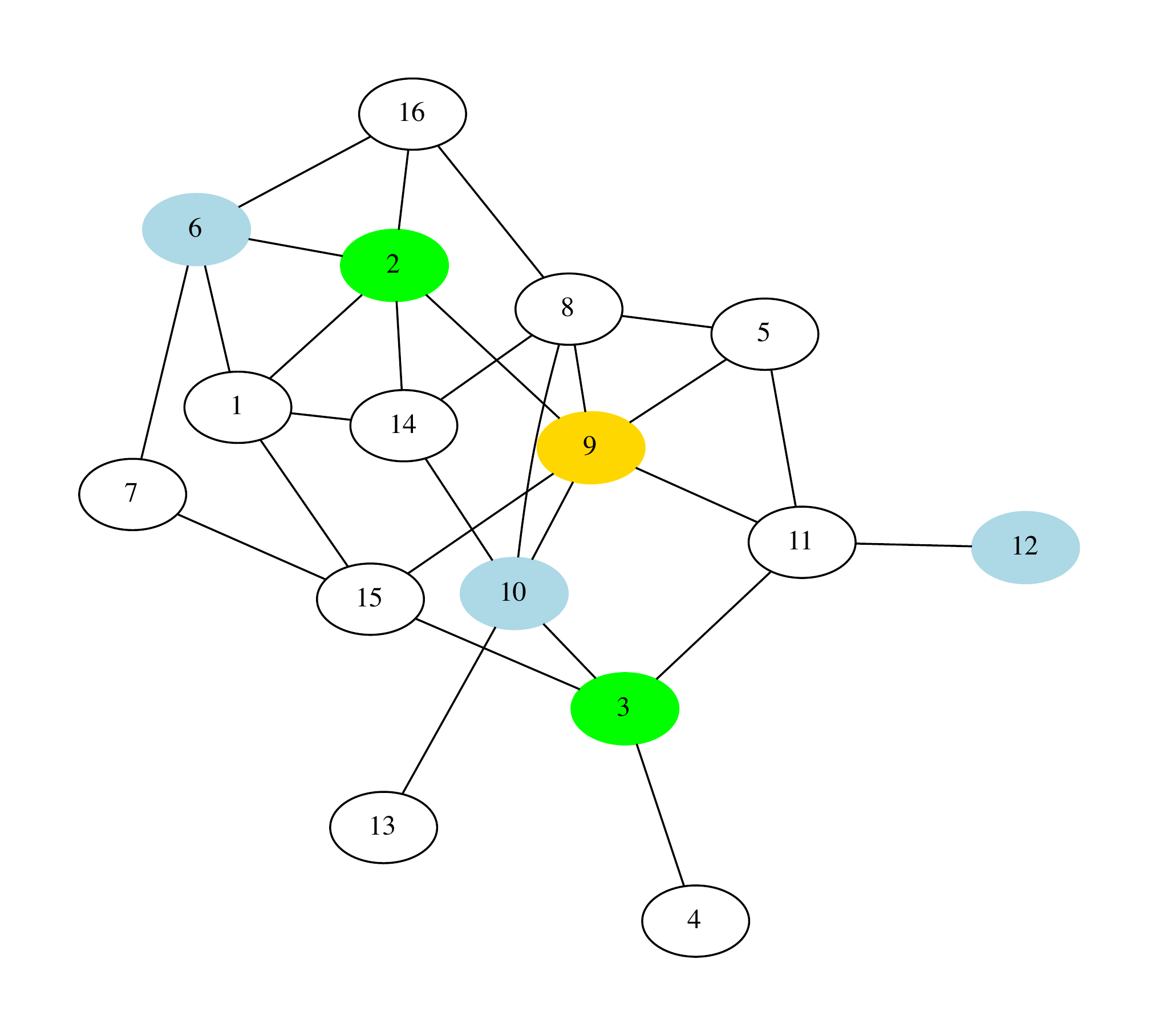}
\caption{Solutions to the optimal leader set for an example graph with sixteen nodes.  For $m=1$ leader, the optimal solution is node 9, shown in yellow.   For $m=2$ leaders, the optimal solution is the set of nodes 2 and 3, shown in green.  For $m=3$ leaders, the optimal solution is the set of nodes 6, 10, and 12, shown in blue.}
\label{fig1}
\end{figure}

To illustrate further, 
we consider the unweighted, undirected, connected graph shown in Figure \ref{fig1}.  The optimal sets of one, two and three leaders are shown in yellow, green and blue, respectively. 
Visually, it is clear that the optimal choice for a single leader (node 9, in yellow) has a central position in the network. In fact, node 9 has the highest information centrality $c_i$ (\ref{infocen}), consistent with Corollary~\ref{cor0} of Section~\ref{sec:interp}, where it is proved that the optimal single leader is the most information central node.  

Interestingly, it is observed that the optimal single leader is not a member of the optimal set of two leaders (nodes 2 and 3, in green). This is due to the fact that the optimal two leaders need to trade off high information centrality as individuals with a joint coverage of the graph (see also Corollary \ref{cor2} in Section~\ref{sec:interp}). For this reason the optimal two leaders are  well connected within the graph and distanced from each other. 

The optimal three leaders (nodes 6, 10, 12, in blue) further illustrate the key trade-off between leaders that are  central and leaders  that cover the graph. 
Although node 12 is not so well connected, its large resistance and biharmonic distances from nodes 6 and 10 make it part of the optimal three-leader set.   That is, the three-node leader set has optimal joint influence on the graph, as encoded by the joint centrality of the set.


The three solutions  illustrate how a leader selection algorithm that first selects a leader and then iteratively adds to the set would result in a sub-optimal leader set for this example and likely in general (see also the example in \cite{Patterson2010}).



\subsection{Joint centrality and two noise-free leaders}
\begin{corollary} \label{cor2} Let $\mathcal{G}$ be an undirected, connected graph of order $n$.   Let  $S_2 = \{s_1,s_2\}$ be a set of two noise-free leaders. Then, the total system error (\ref{arn}) for the system dynamics (\ref{vec1}) is
\begin{align}
\sum_{i=1}^N \Sigma_{ii} = \frac{\sigma^2}{2}\left(\frac{n}{\rho_{S_2}} \right), \label{2nf}
\end{align}
where $\rho_{S_2}$ is the joint centrality of  $S_2$ given by (\ref{rhom}), which specializes to
 \begin{align}
 \rho_{S_2} =n \Big(\frac{K_f}{n} +\frac{nL^+_{s_1,s_1}L^+_{s_2,s_2} - n{L^+_{s_1,s_2}}^2 - \gamma_{s_1,s_2}}{r_{s_1,s_2}} \Big)^{-1}. \label{rho2}
 \end{align}
 The optimal leader set is $S_2^* = \{s_1^*,s_2^*\}= \arg\max_{s_1,s_2} \rho_{S_2}$, the two nodes with the maximal joint centrality.
 \end{corollary}
 
 \proof
 In the case of two leaders, $G = \frac{1}{r_{s_1,s_2}}$.  Equation (\ref{2nf}) follows directly from simplification of  (\ref{first}) and (\ref{rhom}) from Theorem~\ref{thm1}. 
 \endproof
 
 \begin{remark}  Following Remark~\ref{rem2}, we see that in the two-leader case, the term det$(G)$det$(L^+_S) = (L^+_{s_1,s_1}L^+_{s_2,s_2} - {L^+_{s_1,s_2}}^2)/r_{s_1,s_2}$, which is small for large leader information centrality and large resistance distance between leaders.   The term Tr$(Q)$ is proportional to $-\gamma_{s_1,s_1}/r_{s_1,s_2}$.   For this term to be small, the biharmonic distance should be large relative to the resistance distance between leaders.
 \end{remark}

\subsection{Joint centrality and two noise-corrupted leaders}
To address the case of two noise-corrupted leaders, where $k < \infty$, we define a $k$-dependent joint centrality of a set of two nodes.  We  then derive the solution to the optimal leader selection problem for two noise-corrupted leaders by calculating the total system error in terms of the $k$-dependent joint centrality of the two-leader set. 

\begin{theorem} \label{thm2} Let $\mathcal{G}$ be an undirected, connected graph of order $n$.  Let  $S_2 = \{s_1,s_2\}$ be a set of two noise-corrupted leaders ($k < \infty$).  
Define $ \rho_{kS_2}$, the $k$-dependent joint centrality of  $S_2$, as
 \begin{align}
 \rho_{kS_2}& =n\Big(\frac{K_f}{n} + \frac{n[1 + k(L^+_{s_1,s_1}+L^+_{s_2,s_2})]}{k(2+kr_{s_1,s_2})}+ \frac{nk^2(L^+_{s_1,s_1}L^+_{s_2,s_2} - {L^+_{s_1,s_2}}^2) - k^2\gamma_{s_1,s_2}}{k(2+kr_{s_1,s_2})}\Big)^{-1}. \label{rhok2}
 \end{align}
 Then, the total system error (\ref{arn}) for the system dynamics (\ref{vec1}) is
 \begin{align}
\sum_{i=1}^N \Sigma_{ii} = \frac{\sigma^2}{2}\left(\frac{n}{\rho_{kS_2}} \right). \label{2nc}
\end{align}
The optimal leader set is $S_2^* = \{s_1^*,s_2^*\}= \arg\max_{s_1,s_2} \rho_{kS_2}$, the two nodes with the maximal $k$-dependent joint centrality.
\end{theorem}
 
Prior to proving Theorem~\ref{thm2}, we state a lemma from \cite{miller1981} that provides a simplification of the Woodbury formula in the case of a rank one update to a matrix.
\begin{lemma} \label{lemma2} \cite{miller1981} For rank one square matrix $H$ and nonsingular $X$ and $X+H$, $(X+H)^{-1}$ can be written as
\begin{align}
(X+H)^{-1} = X^{-1}-\frac{1}{1+g}X^{-1}HX^{-1}, \label{inv2}
\end{align}
where $g = \mathrm{Tr}(HX^{-1})$.  
\end{lemma}

\proof (Theorem~\ref{thm2}).  Let $K_1$, $K_2$ be rank one matrices with $K_{1_{s_1,s_1}} = k$,  $K_{2_{s_2,s_2}} = k$ where $k>0$ and all other elements of $K_1$, $K_2$ are zero. Let $K = K_1+K_2$ and $N = L+K_1$.  Then, $M = L+K = N+K_2$. 

By applying  Lemma~\ref{lemma2}, we compute
\begin{align}
M^{-1} &= (N+K_2)^{-1} \nonumber\\
& = N^{-1} - \frac{1}{1+\mathrm{Tr}(K_2N^{-1})}N^{-1}K_2N^{-1}. \label{inv1}
\end{align}
By (\ref{nmat}) 
\begin{align}
\mathrm{Tr}(K_2N^{-1}) & = 1+ kL_{s_2,s_2}^+ - 2kL_{s_2,s_1}^++kL_{s_1,s_1}^+  \nonumber \\
& = 1+k\:r_{s_1,s_2}. \label{trace}
\end{align}
Plugging (\ref{trace}) into (\ref{inv1}) yields  total system error (\ref{arn}):
\begin{align}
\sum_{i = 1}^n & \Sigma_{i,i} = \frac{\sigma^2}{2} \sum_{i=1}^n M_{i,i}^{-1} = \frac{\sigma^2}{2} \sum_{i=1}^n \left( N_{i,i}^{-1} - \frac{1}{2 + k r_{s_1,s_2}}(N^{-1}K_2N^{-1})_{i,i} \right).
\end{align}
Expanding $N^{-1}$ in terms of $L^+$ and applying (\ref{Lp2}) and (\ref{Lp3}) gives 
\begin{align}
\sum_{i = 1}^n M_{i,i}^{-1} =& \frac{n}{k}+\frac{K_f}{n} +nL_{s_1,s_1} -\frac{1}{2+k\:r_{s,p}}\Big(k\sum_{i=1}^n(L_{i,s_1}^+ - L_{i,s_2}^+)^2 +nk(L_{s_1,s_2}^+)^2 
- 2nL_{s_1,s_2}^+- \nonumber \\
&2nkL_{s_1,s_1}^+L_{s_1,s_2}^+ +2nL_{s_1,s_1}^++nk(L_{s_1,s_1}^+)^2 +\frac{n}{k} \Big). \label{twotwo}
\end{align}
Rearranging terms and substituting from  (\ref{s1}) results in
\begin{align}
\sum_{i = 1}^n  \Sigma_{i,i} &=   \frac{\sigma^2}{2}\Big(\frac{K_f}{n} + \frac{n + nk(L^+_{s_1,s_1}+L^+_{s_2,s_2})}{k(2+kr_{s_1,s_2})} +\frac{nk^2(L^+_{s_1,s_1}L^+_{s_2,s_2} - {L^+_{s_1,s_2}}^2) - k^2\gamma_{s_1,s_2}}{k(2+kr_{s_1,s_2})})  \Big)\nonumber \\
&= \frac{\sigma^2}{2}\Big(\frac{n}{\rho_{kS_2}}\Big).
\end{align}

\endproof

We observe that the $k$-dependent joint centrality ($\rho_{kS_2}$ from Theorem~\ref{thm2})  plays the same role in determining total system error with noise-corrupted leaders (\ref{2nc}) as joint centrality ($\rho_{S_2}$ from Corollary~\ref{cor2}) plays in determining total system error with noise-free leaders (\ref{2nf}).  

\begin{remark} \textcolor{black}{We note that in the limit as $k \rightarrow \infty$, $\rho_{kS_2}$ approaches $\frac{1}{2}(K_f + n^2(L_{s_1,s_1} + L_{s_2,s_2}))$} Further, as expected, in the limit as $k \rightarrow \infty$ we see that $\rho_{kS_2}$  approaches $\rho_{S_2}$. \textcolor{black}{Therefore one can observe that when the leaders nodes have an infinitely small amount of feedback on their distance from the signal the optimal leader choices are simply the two most information central nodes. As $k$ increases, however, the optimal leader pair changes to satisfy a trade of between information centrality and coverage of the graph.}
\end{remark}

 \subsection{Optimal two noise-corrupted leaders on a cycle}
\begin{corollary} \label{thmcyc} Let $\mathcal{G}$ be an undirected, unweighted cycle graph of order $n$ where $n$ is even.   Let  $S_2 = \{s_1,s_2\}$ be a set of two noise-corrupted leaders ($k < 0$).  The  optimal leader set $S^*$ is any two nodes with maximal resistance distance $r_{s_1,s_2} = \frac{n}{4}$, which corresponds to geodesic distance $d_{s_1,s_2} = \frac{n}{2}$  and antipodal nodes. 
\end{corollary}
 
 \proof See Appendix~\ref{athmcyc}. \endproof
 
 The result is the same for two noise-free leaders on a cycle, i.e., any antipodal pair of leaders solves the optimal two noise-free leader selection problem. \textcolor{black}{In fact, for any network in which all nodes have equivalent information centrality the optimal pair of leaders will be the same for the $k < \infty$ and $k \rightarrow \infty$ cases, and these leaders will be any pair with maximum resistance distance.} Further, we prove in Appendix~\ref{amcyc} that the optimal solution to the general case of $m$ noise-free leaders on a cycle corresponds to the uniform distribution of the leaders around the cycle.    This is consistent with maximizing joint centrality, interpreted as a trade-off between maximizing information centrality of leaders and maximizing coverage of the graph by the leader set:  since every node in a cycle has the same information centrality, the optimal solution only needs to maximize coverage, equivalent to the uniform distribution of leaders around the cycle.
 
 \subsection{Optimal two noise-free leaders on a path} 

\begin{corollary} \label{corpath} Let $\mathcal{G}$ be an undirected, unweighted path graph of order $n$, which is the cycle graph with one link removed.  Let  $S_2 = \{s_1,s_2\}$ be a set of two noise-free leaders.  The optimal leader set $S^*$ is  $s_1^* = \mathrm{rnd}(\frac{n}{5} + \frac{1}{2})$ and $s_2^* = \mathrm{rnd}(\frac{4n}{5} + \frac{1}{2})$, where $\mathrm{rnd}$ is rounding to the closest integer.
\end{corollary}

\proof See Appendix~\ref{acorpath}. \endproof

We observe that for large $n$, the optimal two leader locations on the path approach 0.2 and 0.8 of the path length (starting from one end). This is in contrast with the cycle, where the optimal two leaders  maintain a distance between each other equal to 0.5 of the number of nodes. Considering that the path is simply a cycle with one edge removed, it is interesting to observe that for large $n$, removing an edge from a cycle will cause the fraction of nodes between the  optimal two leaders to increase from 0.5 to 0.6.  That is, the optimal two leaders in the path are more spread out towards the two endpoints.    The locations of the optimal two leaders in the path can be understood to be the optimal solution to the trade-off between high information centrality of  two symmetrically distributed leaders, which increases as the two leaders get closer to midpoint and thus to each other,  and good coverage, which requires the  two leaders to be sufficiently distant from each other.   The optimal two-leader set does not include the optimal single leader set, which is the node at the midpoint of the path, following Corollary~\ref{cor0} of Section~\ref{sec:interp}.

 
 
 
 \subsection{Optimal Selection of a Single Noise-Corrupted and Noise-Free Leader} 
\begin{corollary} \label{cor0} Let $\mathcal{G}$ be an undirected, connected graph of order $n$.  Let  $S = \{s\}$ be a set of one noise-corrupted leader ($k < \infty$)  with information centrality $c_s$. Then, the total system error (\ref{arn}) for the system dynamics (\ref{vec1}) is
\begin{align}
\sum_{i = 1}^n \Sigma_{i,i} = \frac{n\sigma^2}{2} \left(\frac{1}{k} + \frac{1}{c_s} \right). \label{firstagain}
\end{align}
If instead the leader set $S$ is noise-free, then the total system error (\ref{arn}) for the system dynamics (\ref{vec1}) is
\begin{align}
\sum_{i = 1}^n \Sigma_{i,i} = \frac{n\sigma^2}{2} \left(\frac{1}{c_s} \right). \label{toterrfree1}
\end{align}
In both the noise-corrupted and the noise-free cases, the optimal leader set $S^* = \{s^*\}= \arg\max_s c_s$ , the node with maximal information centrality $c_{s^*}$.  
\end{corollary}

\proof
For a single leader we only need to consider a rank-one update to the pseudoinverse of $L$.  From (\ref{nmat}), where $l_1 = s$, this is
\begin{align}
N^{-1} =  L^{-1}  -L_{s}^+\mathbf{1_n}^T -\mathbf{1_n}L_{s}^{+T}  + \frac{ (1+kL_{s,s}^+)}{k}\mathbf{1_n}^T\mathbf{1_n}. \label{nmat2}
\end{align}
Summing the diagonal elements of  (\ref{nmat2})  and applying  (\ref{Lp}), (\ref{Lp2}), (\ref{Lp3}) yields
\begin{align}
\sum_{i=1}^n N_{i,i}^{-1} &= \frac{K_f}{n} + \frac{n}{k} + n\left(\frac{1}{c_s} - \frac{K_f}{n^2}\right) = \frac{n}{k} + \frac{n}{c_s}. \label{minv}
\end{align}
Subsequently substituting into (\ref{arn}) gives the total system error 
\begin{align}
\sum_{i = 1}^n \Sigma_{i,i} = \frac{n\sigma^2}{2} \left(\frac{1}{k} + \frac{1}{c_s} \right). \label{err1}
\end{align}
To get the total system error in the case of one noise-free leader, we take the limit of (\ref{err1}) as $k \rightarrow \infty$, which gives
\begin{align}
\lim_{k\rightarrow \infty} \sum_{i = 1}^n \Sigma_{i,i} = \lim_{k\rightarrow \infty} \frac{n\sigma^2}{2} \left(\frac{1}{k} + \frac{1}{c_s} \right) = \frac{n\sigma^2}{2} \left(\frac{1}{c_s} \right). \label{err2}
\end{align}
The total system error in (\ref{err1}) and in (\ref{err2}) is minimized when the leader has the highest information centrality. \endproof



\subsection{\textcolor{black}{Comparisons to greedy methods}}
\textcolor{black}{As shown in the example above, in general a greedy leader selection method will result in a sub-optimal leader set since the optimal set of $m$ leaders does not necessarily include the optimal set of $m-1$ leaders. However, an advantage of  greedy leader selection methods is that they can be adapted to optimize other quantities, whereas joint centrality was derived specifically to minimize total system error. Conversely, a benefit of calculating joint centrality is that it provides insight into the relevant graph characteristics of the optimal leader set. To demonstrate this take, for example, the cycle graph. We prove in Appendix B that the optimal solution for $m$ noise-free leaders is where the leaders are uniformly distributed and thus this leader set will maximize joint centrality. A greedy method will give the optimal solution for only $m = 2^a$ where $a =0,1,2,3...$ . Because of this, insight into the simple configuration of the optimal leader set is obscured and strictly sub-optimal leader sets are obtained for $m\neq 2^a$ with a greedy approach.}

\subsection{\textcolor{black}{Connections to controllability of networks}}
\textcolor{black}{The problem at hand draws many similarities to the problem of controllability of networks, which has been studied by \cite{summers},  \cite{chapman},  \cite{bullo2014}, among others. When considering network controllability, one has the state equation $\dot{\mathbf{x}} = A\mathbf{x}+B\mathbf{u}$, leading to the calculation of the controllability grammian as 
$W_c = \int_0^\infty e^{A\tau}BB^Te^{A^T \tau} d\tau$. The problem is then to choose $B$ to satisfy a given problem objective. For example this may be to choose $B$ such that $W_c$ is positive definite or achieve a bound the smallest eigenvalue of $W_c$ or to minimize the average control effort given by the trace of $W_c^{-1}$. In contrast, the leader selection problem discussed in this paper considers a fixed $B$ matrix that represents independent noise on each node and seeks to modify the diagonal elements of $A$ such that the trace of $W_c $, equivalently total system error, is minimized. One can see that total system error is in fact a controllability grammian by considering the noise, $\sigma dW$, to be inputs to each node of the system. Therefore total system error is a measure of how robust the system is at the consensus state to random perturbations. The distinction between choosing the $B$ matrix and augmenting diagonal elements of the $A$ matrix is subtle, yet leads to a fundamental difference in both the approach to solving the problem and the solution outcome.
 }

\section{Joint Centrality and Synthetic Lethality in Saccharomyces cerevisiae} \label{sec:special}
To further investigate joint centrality of a set of nodes, we apply it in the analysis of synthetically lethal (SL) genes of the functional gene network of Saccharomyces Cerevisiae, also known as baker's yeast. A functional gene network is one in which nodes in the network represent genes and edges between pairs of nodes represent the function or process by which the pair of genes interact.  S. Cerevisiae has served as a platform for studying genetics of human diseases and is therefore an important model for biological studies \cite{yeastnet}. Here, we focus on instances of synthetic lethality, which occur when the deletion of two genes (A and B) is lethal to the organism and the deletion of A alone or B alone is not lethal. 

\begin{figure}[h!]
\centering
\includegraphics[width=3.5in]{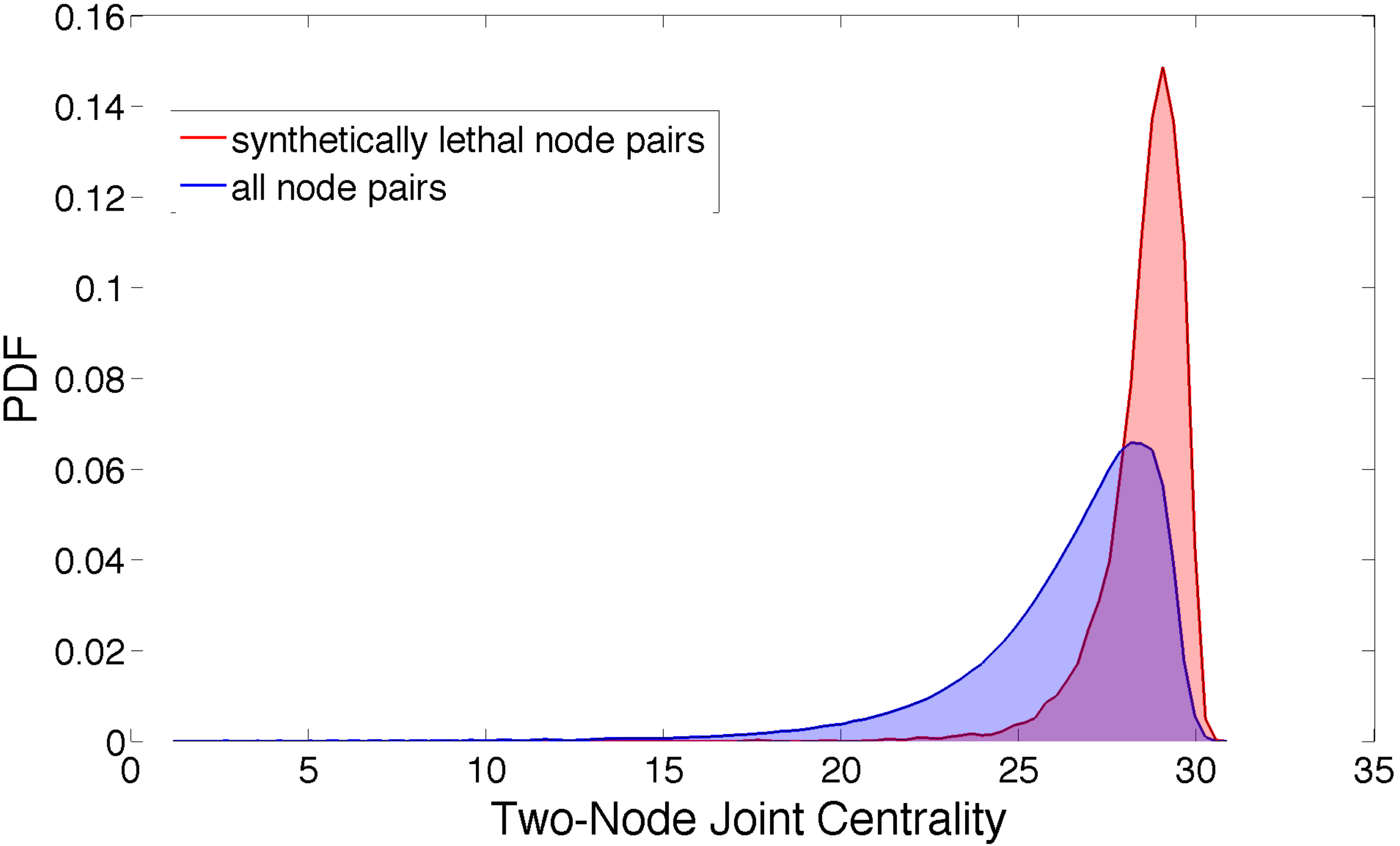} 
\caption{Distribution of two-node joint centrality for every node pair (blue) in the functional gene network of S. Cerevisiae and distribution of two-node joint centrality of synthetically lethal node pairs (red).}
\label{fig2}
\end{figure}

Using the probabilistic functional gene network of S. Cerevisiae from \cite{yeastnet}, (5808 genes with 362,421 edges that represent functional couplings) we calculated the two-node joint centrality for every pair of genes in the network. Then we applied experimental interaction data from the BioGrid database to identify SL pairs of nodes \cite{biogrid}. Figure \ref{fig2} shows the probability distribution function of two-node joint centrality for all pairs of genes (blue) against the probability distribution function of two-node joint centrality for SL  pairs of genes (red). The distributions were constructed by fitting non-parametric distributions with a normal kernel function to normalized histograms of joint centrality calculations for all node pairs and  for all SL node pairs. 

A clear distinction between the two distributions in Figure~\ref{fig2} is apparent. The distribution of two-node joint centralities for SL node pairs is more highly skewed towards high values of joint centrality than the distribution of two-node joint centralities for all node pairs.

We note that SL pairs of nodes are also distinguishable from all other pairs due to their having a higher average degree. This is expected, however, as there is likely a research bias towards testing high degree nodes for synthetic lethality (the set of SL pairs is not necessarily the complete set but rather the set that has been identified thus far).   Accordingly, we do not suggest that joint centrality is the only way to predict possible SL pairs.  Instead, we suggest that two-node joint centrality provides a natural measure for predicting SL pairs, because it takes into account the joint influence of a pair of nodes on the entire network.  In contrast, a measure of pairwise average degree only considers independent, local interactions.

\section{Final Remarks}
\label{sec:final}

In this paper we examine the optimal leader selection problem in a leader-follower network dynamic subject to stochastic disturbances.   The objective is for the network to track an external, unknown signal, where  leaders can take measurements of the external signal but followers must rely only on their measurements of their neighbors. Performance is defined as the inverse of total steady-state error of the system about an external, unknown signal to be tracked, and the optimal set of $m$ leaders maximizes performance over all possible sets of $m$ leaders.   

Our approach is to derive total system error as a function of properties of the underlying network graph.  To do so we define the joint centrality of a set of nodes, such that total system error is proportional to joint centrality.  We prove that the optimal leader set corresponds to the set of $m$ leaders with maximal joint centrality. We show that joint centrality of a set is a function of information centrality of the nodes in the  set and resistance distances and biharmonic distances between pairs of  nodes in the set.  We discuss how the optimal solution is the set of leader nodes that trades off high information centrality with a good coverage of the graph.  

We specialize the results to two noise-corrupted leaders and two noise-free leaders, and we solve for the explicit optimal two-leader solution in the case of the cycle graph and of the path graph.  We also specialize the results to the case of one optimal leader and show that in the noise-corrupted and noise-free cases  the optimal single leader is the node with highest  information centrality. Finally, we provide additional illustration of joint centrality and its more general applicability by using it in the analysis of synthetically lethal gene pairs in a functional gene network.    

Our optimal leader selection results are relevant both to control design, e.g., enabling accuracy and efficiency in sensor networks, and to analysis, e.g., finding conditions that yield the high performance observed in collective animal behavior.  One future direction is to extend the optimal leader selection results of this paper to directed networks by applying the definition in \cite{dres1, dres2} of effective resistance in directed graphs towards a definition of joint centrality in  directed graphs. 
Another compelling future direction is to derive distributed, on-line algorithms that solve the optimal leader selection problem, leveraging our solutions that depend on measures of the graph.   


\appendices
\section{Proof of Corollary~\ref{thmcyc}} \label{athmcyc}
\proof
For a circulant graph, $L^+_{s_1,s_1} = L^+_{s_2,s_2}= L^+_{s,s}$ and thus $\gamma_{s_1,s_2} = \frac{1}{4}\sum_{i=1}^n(r_{i,s_1}-r_{i,s_2})^2$. The $k$-dependent joint centrality $\rho_{kS_2}$ (\ref{rhok2}) simplifies to
\begin{align}
& \rho_{kS_2} = \nonumber \\
 & \;\;n \Big(\frac{K_f}{n} +\frac{n{L^+_{s,s}}^2- n{L^+_{s_1,s_2}}^2 - \sum_{i=1}^n(r_{i,s_1}-r_{i,s_2})^2}{4r_{s_1,s_2}} \Big)^{-1}.
 \end{align}
 By applying  (\ref{res}) and re-arranging terms we have
 \begin{align}
  \rho_{kS_2} &= \frac{n^2}{4}\Big(\frac{K_f}{n^2}+\frac{2}{k} + 4 L^+_{s,s} - r_{s_1,s_2} - \frac{k}{4}\frac{\sum_{i=1}^n(r_{i,s_1}-r_{i,s_2})^2}{2+kr_{s_1,s_2}}\Big)^{-1}. \label{circ1}
  \end{align}
Using the electric circuit analog of resistance distance and applying Kirchhoff's laws, the resistance distance between any two nodes in a cycle can be written as
\begin{align}
\frac{1}{r_{i,j}} = \frac{1}{d_{i,j}} + \frac{1}{n-d_{i,j}},
 \label{rest}
\end{align}
where $d_{i,j}$ is the geodesic distance between nodes $i$ and $j$. The maximum resistance distance is  $r_{i,j} = \frac{n}{4}$, which is obtained between two nodes with $d_{i,j} = \frac{n}{2}$.

Simplifying the $\sum_{i=1}^n(r_{i,s_1}^+ - r_{i,s_2}^+)^2$ term of  (\ref{circ1}) by inserting (\ref{rest}) gives
\begin{align}
\sum_{i=1}^n(r_{i,s_1} &- r_{i,s_2})^2 =\sum_{i=1}^n \Big (d_{i,s_1}-d_{i,s_2}+ \frac{d_{i,s_2}^2-d_{i,s_1}^2}{n}\Big)^2 \nonumber \\
 & = \frac{d_{s_1,s_2}(d_{s_1,s_2}-n)(d_{s_1,s_2}^2-nd_{s_1,s_2}-2)}{3n}. \label{sumsq2}
\end{align}
Substituting (\ref{sumsq2}) into (\ref{circ1}) results in
\begin{align}
\rho_{kS_2} &=  \frac{n^2}{4}\Big(\frac{K_f}{n^2}+\frac{2}{k} + 4 L^+_{s,s} - \frac{d_{s_1,s_2}(n-d_{s_1,s_2})}{n} - \frac{k d_{s_1,s_2}(d_{s_1,s_2}-n)(d_{s_1,s_2}^2-nd_{s_1,s_2}-2)}{6n (2n +k d_{s_1,s_2}(n-d_{s_1,s_2}))}\Big)^{-1}.  \label{ninv5}
\end{align}

To determine how $\rho_{kS_2}$ changes as a function of $d_{s_1,s_2}$, we take the partial derivative of (\ref{ninv5}) with respect to $d_{s_1,s_2}$ to give
\begin{align}
\frac{\partial \rho_{kS_2}^{-1}}{\partial d_{s_1,s_2}}=& - \frac{1}{4}(n-2d_{s_1,s_2}) -\frac{nk [2(- d_{s_1,s_2} +  d_{s_1,s_2}^3)  + (1 - 3 d_{s_1,s_2}^2) n +  d_{s_1,s_2} n^2] }{3 (2 n + d_{s_1,s_2} k (-d_{s_1,s_2} + n))^2} \nonumber \\
& -  \frac{k^2 [-2d_{s_1,s_2}^5 + 5 d_{s_1,s_2}^4 n - 4 d_{s_1,s_2}^3 n^2 + d_{s_1,s_2}^2 n^3]}{12 (2 n + d_{s_1,s_2} k (-d_{s_1,s_2} + n))^2}. \label{deriv1}
 \end{align}
Since $d_{s_1,s_2} \leq \frac{n}{2}$, the first term of (\ref{deriv1}) will always be nonpositive. Additionally, it can be shown algebraically that  for $n>3$ the two bracketed expressions in the second and third terms will be greater than zero. Therefore $\rho_{kS_2}^{-1}$ decreases as $d_{s_1,s_2}$ increases, reaching its minimum at the maximal value of $d_{s_1,s_2} = \frac{n}{2}$, corresponding to $r_{s_1,s_2} = \frac{n}{4}$.  
 \endproof
 
 \section{Optimal $m$ noise-free leaders on a cycle} \label{amcyc}
In the case where the network graph is a cycle, we can use the cyclic structure of the graph Laplacian  to explicitly solve for the optimal locations of  $m$ noise-free leaders.
\begin{theorem} Let $\mathcal{G}$ be an undirected, unweighted cycle graph of order $n$.   Let $m < n$  such that $p = n/m$ is an integer.  Let  $S$ be a set of $m$ noise-free leaders. Then, an optimal leader set $S^*$ is any set $S$ where the leaders are uniformly distributed around the cycle, i.e., the geodesic distance between any leader and each of the other two closest leaders is  $d_{s_a,s_b} = p$.
\end{theorem}

\proof
 We begin by assuming $m$ nodes on the cycle have been selected as leaders and let $M = L+K$ where $K$ is a matrix with a value of $k$ in the entries along the main diagonal corresponding to the leader nodes and zeros elsewhere. We partition $M$ in the usual way. Since we are assuming noise-free leaders, to compute total system error we need only to consider the sum of the diagonal elements of the inverse of the submatrix $M_{F}$. $M_{F}$ can be written as a block diagonal matrix where each block corresponds to a set of connected follower nodes between two leader nodes. Each block, $M_{{F}_i}$ will itself be a tridiagonal matrix of the form
\begin{align}
M_{{F}_i} =  \begin{bmatrix} 2& -1 &  & 0  \\ -1 & 2 & & &\\ 0 & & \ddots & -1\\ 0 & & -1 & 2 \end{bmatrix}.
\end{align}
In the case where there is one follower node in between two leader nodes the corresponding diagonal block in $M_{F}$ will be one element with an entry of $2$. 

Similar to previous sections, total system error for noise-free leaders  will be proportional to the trace of $M_{F}^{-1}$, which here is equivalent to the total sum of eigenvalues of each $M_{{F}_i}^{-1}$. 
By \cite{tritoep} we have that the eigenvalues of  $M_{F_i}^{-1}$ are 
\begin{align}
\lambda_{z_{ij}} = \frac{1}{2 - 2\cos \left(j \frac{\pi}{w_i+1}\right)} \; \; \;\; j = 1,...., w_i
\end{align}
where $w_i$ is dimension of $M_{{F}_i}$. The average value of the eigenvalues of a block is then
\begin{align}
\bar\lambda_{z_i} = \sum_{j=1}^{w_i} \lambda_{z_{ij}}= \frac{1}{6}w_i +\frac{1}{3}.
\end{align}

Therefore, minimizing the total sum of eigenvalues is equivalent to minimizing the sum over $i$ of $w_i^2$. It follows that the minimum is achieved when $w_1 = w_2 = w_3 = ...$, or in other words when the dimension of each block is the same.  This corresponds to the leaders being evenly distributed around the cycle with shortest distances between leaders equal to $d_{s_1,s_2} = \frac{n}{2}$. 
\endproof

\section{Proof of corollary \ref{corpath}} \label{acorpath}
\proof
Resistance distance in a path graph simplifies to $ r_{i,j} = \|i-j\| $ and
\begin{align} 
L_{j,j} &= \frac{\sum_{i=1}^n r_{i,j}}{n} -\frac{K_f}{n^2}\nonumber \\
&= \frac{(n-j)(1+n-j) - j +j^2}{2n} -\frac{K_f}{n^2}. \label{p2}
\end{align}
Substitution of (\ref{p2}) into the expression (\ref{rho2}) for $\rho_{S_2}$, where without loss of generality we take  $s_2 >s_1$, gives
\begin{align}
\rho_{S_2}^{-1} &= \frac{1}{n}\Big(-\frac{1}{6}+\frac{n +n^2 -s_1 -s_2}{4} +\frac{(2s_1^2+2s_2^2-s_2(3n+s_1))}{3}\Big).  \label{path2}
\end{align}

We then take partial derivatives of (\ref{path2})  with respect to $s_1$ and $s_2$ to find  the minimum of $\rho_{S_2}^{-1}$ to be $s_1 =  \text{rnd}(\frac{n}{5} + \frac{1}{2})$ and $s_2 = \text{rnd}(\frac{4n}{5} + \frac{1}{2})$.   The rounding of $s_1$ and $s_2$ can be checked by observing from  (\ref{path2}) that the level sets of $\rho_{S_2}^{-1}$ are ellipses in $s_1$, $s_2$.  Computing the semi-axis lengths of the ellipses shows that the nearest integer values of $s_1$ and $s_2$ that minimize $\rho_{S_2}^{-1}$ indeed determine the optimal leader set.   
\endproof

\FloatBarrier
\bibliographystyle{IEEEtran} 
\bibliography{LeaderConsensus2}

\end{document}